\documentstyle[leqno]{article}
\begin{document}
\renewcommand{\thesection}{\arabic{section}.}
\renewcommand{\theequation}{\arabic{section}.\arabic{equation}}
\newcommand{\be}{\begin{eqnarray}}
\newcommand{\en}{\end{eqnarray}}
\newcommand{\no}{\nonumber}
\newcommand{\la}{\lambda}
\newcommand{\laa}{\Lambda}
\newcommand{\ep}{\epsilon}
\newcommand{\de}{\delta}
\newcommand{\D}{\Delta}
\newcommand{\pl}{\parallel}
\newcommand{\ov}{\overline}
\newcommand{\bet}{\beta}
\newcommand{\al}{\alpha}
\newcommand{\fr}{\frac}
\newcommand{\pa}{\partial}
\newcommand{\we}{\wedge}
\newcommand{\om}{\Omega}
\newcommand{\na}{\nabla}
\newcommand{\lan}{\langle}
\newcommand{\ra}{\rangle}
\newcommand{\vs}{\vskip0.3cm}
\renewcommand{\thefootnote}{}
\title {Universal Inequalities for Eigenvalues of the Buckling Problem of Arbitrary Order }
\footnotetext{2000 {\it Mathematics Subject Classification }: 35P15, 53C20, 53C42, 58G25
\hspace*{2ex}Key words and phrases: Universal bounds,  eigenvalues, buckling problem,  Euclidean space, sphere.}
\author{J\"urgen Jost, Xianqing Li-Jost, Qiaoling Wang, Changyu Xia }\date{} \maketitle
\begin{abstract}
We investigate the eigenvalues of the buckling problem of arbitrary
order on compact domains  in  Euclidean spaces and spheres. We obtain
universal bounds for the $k$th eigenvalue in terms of the lower  eigenvalues independently  of the particular geometry of the domain.
\end{abstract}


\section{Introduction}
Let $\om$ be a connected bounded domain with smooth boundary in an  $n(\geq 2)$-dimensional Euclidean space ${I\!\!R}^n$ and let $\nu$ be the outward unit normal vector field of $\pa \om$.
Denote by $\Delta $ the Laplacian operator on ${I\!\!R}^n.$ Let us consider the following well-known eigenvalue problems :
\be& &
\Delta u=-\la u \ \ \ \ \ {\rm in \ \ } \om,  \ \ \ \
u=0, \ \ \ \ \ \ \ \ \ \ \ \ {\rm on \ \ \  } \pa \om,
\\  & &
\Delta^2 u= \eta u \ \ \ \ \  \ {\rm in \ \ } \om,  \ \ \ \
u=\fr{\pa u}{\pa \nu }=0, \ \ \ \  {\rm on \ \ \ } \pa \om,
\\  & &
\Delta^2 u= -\laa \Delta u \ \ \ \ {\rm in \ \ } \om,  \ \ \ \
u=\fr{\pa u}{\pa \nu }=0, \ \ \ \  {\rm on \ \ } \pa \om.
\en
They are called the {\it  fixed membrane problem}; the {\it clamped plate problem} and
the {\it bucking problem}, respectively.
Let
\be\no
& & 0<\la_1<\la_2\leq\la_3\leq\cdots,\\
\no & & 0<\eta_1\leq\eta_2\leq\la_3\leq\cdots,\\
\no & & 0<\laa_1\leq\laa_2\leq\laa_3\leq\cdots
\en
denote the successive eigenvalues for (1.1), (1.2) and (1.3), respectively.
Here each eigenvalue is repeated according to its
multiplicity. Deriving bounds for these (and other) eigenvalues is an
important theme of mathematical analysis. In most cases, eigenvalues
are controlled by the geometry of the underlying domain, the
$n$-dimensional ball often representing an extremal case. On the other
hand, it has been found that one can also control higher eigenvalues
in terms of lower ones, completely independently of the geometry of
the domain (apart from its dimension). Such eigenvalue bounds are
called universal.
Universal bounds for the eigenvalues $\la_{k+1}$, $\eta_{k+1}$ and
$\laa_{k+1}$ have been derived by many mathematicians, and we shall
now recall the pertinent results.
Payne, P\'olya and Weinberger ([PPW1], [PPW2]) proved the bound
\be
\la_{k+1}-\la_k\leq \fr 2k\sum_{i=1}^k\la_i, \ \ k=1, 2,\cdots,
\en
for $\om\subset {I\!\!R}^2$. This result easily extends to $\om\subset {I\!\!R}^n$ as
\be
\la_{k+1}-\la_k\leq \fr 4{kn}\sum_{i=1}^k\la_i, \ \ k=1, 2,\cdots,
\en
In 1980, Hile and Protter [HP] proved
\be
\sum_{i=1}^k\fr{\la_i}{\la_{k+1}-\la_i}\geq \fr{kn}4, \ \ \ {\rm for} \ k=1,2,\cdots.
\en
In 1991, Yang proved the following much stronger inequality [Y]:
\be
\sum_{i=1}^k(\la_{k+1}-\la_i)\left(\la_{k+1}-\left(1+\fr 4n\right)\la_i\right)\leq 0, \ \ \
{\rm for \ } k=1, 2,\cdots.
\en
The inequality (1.7)  is the strongest of the classical inequalities that are derived following the scheme devised by Payne-P\'olya-Weinberger. Yang's inequality provided a marked improvement for eigenvalues of large index.
It should be also  mentioned that the development of Yang's inequality
came to fruition only thanks to the work of M. S. Ashbaugh [A2] and
that of Harrell-Stubbe [HS]. In fact, it was Harrell-Stubbe who first
explained the key commutator facts behind the ``trick'' introduced by
H. C. Yang in the traditional Payne-P\'olya-Weinberger scheme and
introduced the Yang inequality to the mathematical physics and
geometry community. This trick was explained in further work of
Ashbaugh (and later in the work of Ashbaugh-Hermi [AH1], [AH2]) as an
instance of the use of the ``optimal Cauchy-Schwarz'' inequality. It
was Ashbaugh who dubbed it the ``Yang inequality''. The optimal
Cauchy-Schwarz trick is what enabled  Cheng-Yang [CY2] and  Wang-Xia [WX1] to extend the earlier work of H. C. Yang to the case of the clamped plate problem for bounded domains of Euclidean space and of minimal submanifolds of the same space, respectively. This is the trick that makes all extensions \`a la H. C. Yang.  The arguments around this trick were later generalized by Harrell [H], Harrell-Michel [HM1], [HM2] and Levitin-Parnovski [LP], following the commutator method via Rayleigh-Ritz.

Consider now the problem (1.3) which is used to describe the critical buckling load of a clamped plate subjected to a uniform compressive force around its boundary.
In 1956,   Payne, P\'olya and Weinberger  proposed in [PPW2] the following
\vskip0.3cm
\noindent
{\bf Problem 1.} {\it  Can one obtain a universal inequality for the
  eigenvalues of the buckling problem (1.3) that is similar to the universal inequalities for the eigenvalues of the fixed membrane problem (1.1)  ?
}
\vskip0.3cm
Ashbaugh [A1] mentioned this problem again. With respect to the above problem,
Payne, P\'olya and Weinberger proved
$$\laa_2/\laa_1<3 \ \ \ \ {\rm for }\ \om\subset I\!\!R^2.$$
For $\om\subset I\!\!R^n$ this reads
$$\laa_2/\laa_1<1+4/n.$$
Subsequently Hile and Yeh [HY] reconsidered this problem obtaining the improved bound
$$
\fr{\laa_2}{\laa_1}\leq \fr{n^2+8n+20}{(n+2)^2}\ \ \ \ \ {\rm for \ } \om\subset I\!\!R^n.
$$
Ashbaugh [A1] proved :
\be
\sum_{i=1}^n\laa_{i+1}\leq(n+4)\laa_1.
\en
Recently, Cheng and Yang introduced a new method to construct trial functions for the problem (1.3) and obtained the following universal inequality [CY3]:
\be
\sum_{i=1}^k(\laa_{k+1}-\laa_i)^2\leq \fr{4(n+2)}{n^2}\sum_{i=1}^k(\laa_{k+1}-\laa_i)\laa_i.
\en
It has been proved in [WX2] that for the problem (1.3) if $\om$ is a bounded connected domain in an $n$-dimensional unit sphere, then the following inequality holds
\be & &
2\sum_{i=1}^k(\laa_{k+1}-\laa_i)^2\\ \no &\leq&
\sum_{i=1}^k(\laa_{k+1}-\laa_i)^2\left(\delta \laa_i+\fr{\delta^2(\laa_i-(n-2))}{4(\delta\laa_i+n-2)}\right)  +\fr 1{\delta}\sum_{i=1}^k(\laa_{k+1}-\laa_i)\left(\laa_i+\fr{(n-2)^2}4\right),
\en
where $\delta $ is any positive constant.

In this paper, we will investigate the eigenvalues of the  buckling problem of higher order:
\be
& & (-\Delta)^l u= -\laa\D u\ \ \ {\rm in} \ \ \om, \\
\no  & &  u|_{\pa \om}=\left. \fr{\pa u}{\pa \nu}\right|_{\pa \om}=\cdots =\left. \fr{\pa^{l-1} u}{\pa \nu^{l-1}}\right|_{\pa \om}=0,
\en
where $ \om$ is a connected bounded domain in a Euclidean space or a unit sphere and $l$ is any integer no less than $2$.

For the eigenvalues of the problem (1.11), Chen-Qian([CQ]) obtained some
upper bounds on the $k$th eigenvalue  in terms of the lower ones
when $k$ is small and $\om$ is contained in a Euclidean space. To the authors' knowledge, there are no universal
inequalities on $\laa_{k}$ in terms of $\laa_1, \cdots, \laa_{k-1}$
for general $k$. The purpose of this paper is to prove such
inequalities. Namely, we will prove \vskip0.3cm {\bf Theorem 1.1.}
{\it Let $l\geq 2$ and  let $\laa_i$ be the $i$-th eigenvalue of the
following eigenvalue problem: \be
& & (-\Delta)^l u= -\laa\D u\ \ \ {\rm in} \ \ \om, \\
\no  & &  u|_{\pa \om}=\left. \fr{\pa u}{\pa \nu}\right|_{\pa \om}=\cdots =\left. \fr{\pa^{l-1} u}{\pa \nu^{l-1}}\right|_{\pa \om}=0.
\en
where $\om$ is a connected bounded domain in an $n$-dimensional Euclidean space   with smooth boundary $\pa \om$ and $\nu$ is the unit outward normal vector field of $\pa \om$. Then for $k=1, \cdots,$
we have
\be
& &\sum_{i=1}^k(\laa_{k+1}-\laa_i)^2\\ \no
&\leq& \fr{ {2(2l^2+(n-4)l+2-n)}^{1/2}}n\left\{ \sum_{i=1}^k  (\laa_{k+1}-\laa_i)^2 \laa_i^{(l-2)/(l-1)}\right\}^{1/2}
\left\{ \sum_{i=1}^k  (\laa_{k+1}-\laa_i)\laa_i^{1/(l-1)}\right\}^{1/2}.
\en
}
\vskip0.3cm
{\bf Remark.} If we take $l=2$ in Theorem 1.1, then we obtain Cheng-Yang's inequality (1.9).
\vskip0.3cm

From Theorem 1, we can obtain  more explicit inequalities which are weaker than (1.13):
\vskip0.3cm
{\bf Corollary 1.1.} {\it Under the same assumptions as in Theorem 1, we have
\be
\laa_{k+1}
&\leq &\fr 1k\sum_{i=1}^k\laa_i+\fr{2(2l^2+(n-4)l+2-n) }{k^2n^2}\left(\sum_{i=1}^k\laa_i^{(l-2)/(l-1)}
\right)\left(\sum_{i=1}^k\laa_i^{1/(l-1)}\right)\\ \no
& & \ \ \ \ +\left\{\left(\fr{2(2l^2+(n-4)l+2-n)}{k^2n^2}\right)^2\left(\sum_{i=1}^k\laa_i^{(l-2)/(l-1)}
\right)^2\left(\sum_{i=1}^k\laa_i^{1/(l-1)}\right)^2\right.\\ \no & & \left.\ \ \ \ \ \ \ \ \ \ -\fr 1k\sum_{i=1}^k\left(\laa_i-\fr 1k\sum_{j=1}^k \la_j\right)^2\right\}^{\fr 12}.
\en
and
\be
\laa_{k+1}
&\leq & \left(1+\fr{2(2l^2+(n-4)l+2-n)}{n^2}\right)\fr 1k\sum_{i=1}^k\laa_i\\
\no & & \ \ \ \ +
\left\{\left(\fr{2(2l^2+(n-4)l+2-n)}{n^2}\fr 1k\sum_{i=1}^k\laa_i\right)^2\right. \\ \no & &\left. \ \ \ \ \ \ \ \ \ \ -\left(1+\fr{4(2l^2+(n-4)l+2-n)}{n^2}\right)\fr 1k\sum_{i=1}^k\left(\laa_i-\fr 1k\sum_{j=1}^k\laa_j\right)^2\right\}^{1/2}
\en
}
\vskip0.3cm
We then prove the following universal inequalities for eigenvalues of the buckling problem of higher orders on spherical domains.
\vskip0.3cm
{\bf Theorem 1.2.} {\it Let $l\geq 2$ and  let $\laa_i$ be the $i$-th eigenvalue of the following eigenvalue problem:
\be
& & (-\Delta)^l u= -\laa\D u\ \ \ {\rm in} \ \ \om, \\
\no  & &  u|_{\pa \om}=\left. \fr{\pa u}{\pa \nu}\right|_{\pa \om}=\cdots =\left. \fr{\pa^{l-1} u}{\pa \nu^{l-1}}\right|_{\pa \om}=0.
\en
where $\om$ is a connected bounded domain in an $n$-dimensional Euclidean sphere   with smooth boundary $\pa \om$ and $\nu$ is the unit outward normal vector field of $\pa \om$. For each $q=0, 1, \cdots, $ define the polynomials $F_q$ and $G_q$ inductively by
\be F_0=G_0=1, \ F_1(t)=t-(n+2), \ G_1(t)=3t+n-2,\en
\be& &
F_q(t)=(2t-2)F_{q-1}(t)-(t^2+2t-n(n-2))F_{q-2}(t), \\  & & G_q(t)=(2t-2)G_{q-1}(t)-(t^2+2t-n(n-2))G_{q-2}(t), \ \ q=2,\cdots
\en
Set
\be
tF_{l-2}(t)- G_{l-2}(t)= t^{l-1}+a_{l-2}t^{l-2}+\cdots +a_{1}t+ a_0.
\en
Let $\delta$ be any positive number and $k$ be a positive integer.
Then we have
\be  & &
\sum_{i=1}^k(\laa_{k+1}-\laa_i)^2\left(2+\fr{n-2}{\laa_i^{1/(l-1)}-(n-2)}\right)
\\ \no &\leq &
 \delta \sum_{i=1}^k(\laa_{k+1}-\laa_i)^2H_i+
\fr 1{\delta} \sum_{i=1}^k(\laa_{k+1}-\laa_i)\left(\laa_i^{1/(l-1)}+\fr{(n-2)^2}4\right),
\en
where
\be
H_i&=&\laa_i^{1/(l-1)} \left( 1
-\fr 1{\laa_i^{1/(l-1)}-(n-2)}\right)+\sum_{j=0}^{l-2}|a_j|\laa_i^{j/(l-1)}
\en
\vskip0.3cm
}
\vskip0.3cm
{\bf Remark.} When $l=2$, it is easy to see that
$$
H_i=1+\laa_i \left(1-\fr 1{\laa_i-(n-2)}\right)
$$
and so the inequality (1.21) in this case  can be written as
\be
& & 2\sum_{i=1}^k(\laa_{k+1}-\laa_i)^2\\ \no
&\leq &
 \sum_{i=1}^k(\laa_{k+1}-\laa_i)^2  \left(\delta+ \delta\laa_i\left(1-\fr 1{\laa_i-(n-2)}\right)-\fr {n-2}{\laa_i-(n-2)}\right)\\ \no & & +
\fr 1{\delta} \sum_{i=1}^k(\laa_{k+1}-\laa_i)\left(\laa_i^{1/(l-1)}+\fr{(n-2)^2}4\right),
\en
Observe that (1.23) is sharper than (1.10) since for any $\delta >0$, we have
$$
\delta+ \delta\laa_i\left(1-\fr 1{\laa_i-(n-2)}\right)-\fr {n-2}{\laa_i-(n-2)}\leq \delta \laa_i+\fr{\delta^2(\laa_i-(n-2))}{4(\delta\laa_i+n-2)}.
$$

\vskip0.3cm
From Theorem 1.2, we can obtain an explicit upper bound on $\laa_{k+1}$ in terms of $\laa_1,\cdots, \laa_k$ which is weaker than (1.21).
\vskip0.3cm
{\bf Corollary 1.2.} {\it Let the  assumptions and the notations be as in Theorem 1.2. It holds

\be
\laa_{k+1}\leq A_{k+1}+\sqrt{A_{k+1}^2-B_{k+1}},
\en
where
\be
& & A_{k+1}=\fr 1k\sum_{i=1}^k\laa_i +\fr 2{kS_k^2}\sum_{i=1}^k T_i, \ B_{k+1}=\fr 1k\sum_{i=1}^k\laa_i^2 +\fr 4{kS_k^2}\sum_{i=1}^k T_i\laa_i, \\
& &
S_k=2+\fr{n-2}{\laa_k^{1/(l-1)}-(n-2)},\ \ T_i= H_i\left(\laa_i^{1/(l-1)}+\fr{(n-2)^2}4\right).
\en
}

\markboth{\hfill J\"urgen Jost, Xianqing Li-Jost, Qiaoling Wang, Changyu Xia \hfill}
{\hfill  Universal inequalities  for Eigenvalues of the Buckling Problem of Any Order
 \hfill}
\section{Proofs of the Results }
\setcounter{equation}{0}
Before proving our results, let us recall a method of constructing trial functions developed by Cheng-Yang (Cf. [CY3], [WX2]). We will state it in a quite general form since we believe that it  could be useful for studying eigenvalues of the buckling problem of high orders on compact domains
of complete submanifolds in a Euclidean space. Let $M$ be an
  $n$-dimensional complete submanifold in an $m$-dimensional Euclidean space $I\!\!R^{m}$.
Denote by $\langle , \rangle$ the canonical metric on $I\!\!R^{m}$ as well as that induced on $M$. Denote by  $\Delta $ and $\nabla$    the Laplacian and the gradient operator of $M$, respectively. Let $\om$ be a bounded connected domain of $M$ with smooth boundary $\pa \om$ and let $\nu $ be the outward unit normal vector field of $\pa\om$.
For functions $f$ and $g$ on $\om$, the {\it Dirichlet inner product $(f, g)_D$} of $f$ and $g$ is given by
\be \no (f, g)_D=\int_{\om}\langle\nabla f, \ \nabla g\rangle.
\en
The Dirichlet norm of a function $f$ is defined by
\be\no
||f||_D=\{(f, f)_D\}^{1/2}=\left(\int_{\om}|\nabla f|^2\right)^{1/2}.
\en
Consider the
eigenvalue problem
\be
& & (-\Delta)^l u=-\la \Delta u\ \ \ {\rm in} \ \ \om, \\
\no  & &  u|_{\pa \om}=\left. \fr{\pa u}{\pa \nu}\right|_{\pa \om}=\cdots =\left. \fr{\pa^{l-1} u}{\pa \nu^{l-1}}\right|_{\pa \om}=0.
\en
Let
\be\no
 0<\laa_1\leq \laa_2\leq\laa_3\leq\cdots,
\en denote the successive eigenvalues, where
each eigenvalue is repeated according to its multiplicity.

Let $u_i$ be the $i$-th orthonormal eigenfunction of the problem (2.1) corresponding to the eigenvalue $\laa_i$, $i=1, 2, \cdots,$ that is, $u_i$ satisfies
\be& &
(-\Delta)^lu_i=-\laa_i \Delta u_i \ \ {\rm in \ \ } \om,  \\ \no
& &
\left. u_i\right|_{\pa \om}=\left. \fr{\pa u_i}{\pa \nu}\right|_{\pa \om}=\cdots\left. \fr{\pa^{l-1} u}{\pa \nu^{l-1}}\right|_{\pa \om}=0,\\  \no
& &  (u_i, u_j)_D =\int_{\om}\langle \nabla u_i, \nabla u_j\rangle=\delta_{ij},\ \ \ \forall\ i, j.
\en
For $k=1, \cdots, l$, let $\nabla^k$ denote the $k$-th covariant
derivative operator on $M$, defined in the usual weak sense via an
integration by parts formula.
 For a function $f$ on $\om$, the squared norm of  $\nabla^k f$ is defined as (cf. [He])
\be
\left|\nabla^kf\right|^2=\sum_{i_1,\cdots, i_k=1}^n\left(\nabla^kf(e_{i_1},\cdots, e_{i_k})\right)^2,
\en
where $e_1,\cdots, e_n$ are orthonormal vector fieds locally defined on $\om$.
Define the Sobolev space $H_l^2(\om)$ by
$$H_l^2(\om)=\{ f:\ f, \ |\nabla f|,\cdots, \left|\nabla^l f\right|\in L^2(\om)\}.
$$
Then $H_l^2(\om)$ is a Hilbert space with respect to the norm $||\cdot||_{l, 2}$:
\be
||f||_{l, 2}=\left(\int_{\om}\left(\sum_{k=0}^l|\nabla^k f|^2\right)\right)^{1/2}.
\en
Consider the subspace $H_{l,D}^2(\om)$ of $H_l^2(\om)$ defined by
$$H_{l,D}^2(\om)=\left\{f\in H_l^2(\om): \ f|_{\pa \om}=\left. \fr{\pa f}{\pa \nu}\right|_{\pa \om}=\cdots\left. \fr{\pa^{l-1} u}{\pa \nu^{l-1}}\right|_{\pa \om}=0\right\}.
$$
The operator $(-\Delta)^l$ defines a self-adjoint operator acting on $H_{l,D}^2(\om)$
with discrete eigenvalues $0<\laa_1\leq\cdots\leq \laa_k\leq\cdots$ for the buckling problem (2.1) and the eigenfunctions $\{u_i\}_{i=1}^{\infty}$ defined in (2.2)
form a complete orthonormal basis for the Hilbert space $H_{l,D}^2(\om)$. If $\phi\in H_{l,D}^2(\om)$ satisfies $(\phi , u_j)_D=0, \ \forall  j=1, 2, \cdots, k$, then the Rayleigh-Ritz inequality tells us that
\be
\laa_{k+1}||\phi ||_D^2\leq \int_{\om} \phi(-\Delta)^l\phi.
\en
For vector-valued functions $F=(f_1, f_2, \cdots, f_{m}), \ G=(g_1, g_2, \cdots, g_{m}):
\om\rightarrow I\!\!R^{m}$, we define an inner product $(F, G)$ by
$$(F, G)\equiv \int_{\om} \langle F, G\rangle =\int_{\om} \sum_{\alpha =1}^{m} f_{\alpha}g_{\alpha}.$$
The norm of $F$ is given by
$$||F||=(F, F)^{1/2}=\left\{\int_{\om}\sum_{\alpha=1}^{m}f_{\alpha}^2\right\}^{1/2}.$$
Let ${\bf H}_1^2(\om)$ be the Hilbert space of vector-valued functions given by
$${\bf H}_1^2(\om)=\left\{ F=(f_1,\cdots, f_{m}): \om\rightarrow I\!\! R^{m};\  f_{\alpha}, \ |\nabla f_{\alpha}|\in L^2(\om), \ {\rm for} \ \alpha=1,\cdots, m\right\}
$$
with norm
$$||F||_1=\left(||F||^2+\int_{\om}\sum_{\alpha=1}^{m}|\nabla f_{\alpha}|^2\right)^{1/2}.
$$
Observe that a vector field on $\om$ can be regarded as a vector-valued function from $\om$ to $I\!\!R^{m}$. Let ${\bf H}_{1, D}^2(\om)\subset  {\bf H}_{l}^2(\om)$ be a subspace of ${\bf H}_l^2(\om)$
spanned by the vector-valued functions $\{ \nabla u_i\}_{i=1}^{\infty}$, which form a complete orthonormal basis of ${\bf H}_{1, D}^2(\om)$. For any $f\in H_{l,D}^2(\om), $ we have $\nabla f\in {\bf H}_{1, D}^2(\om)$ and for any $X\in {\bf H}_{1, D}^2(\om)$, there exists a function $f\in H_{l,D}^2(\om)$ such that $X=\nabla f$.

\vskip0.3cm
{\bf Lemma 2.1.} {\it Let $u_i$ and $\laa_i, i=1, 2, \cdots , $ be as in (2.2), then
\be
0\leq \int_M u_i(-\Delta )^k u_i\leq \laa_i^{(k-1)/(l-1)}, \ \ k=1,\cdots, l-1.
\en
}

{\it Proof of Lemma 2.1.} When $k\in \{1,\cdots, l-1\}$ is even, we have
\be
\int_M u_i(-\D)^k u_i =\int_M u_i\D^k u_i=\int_M \left(\D^{k/2}u_i\right)^2\geq 0.
\en
On the other hand, if $k\in \{1,\cdots, l-1\}$ is odd, it holds
\be \no
\int_M u_i(-\D)^k u_i &=& -\int_M u_i\D^k u_i
\\ \no
&=& -\int_M \D^{(k-1)/2}u_i\D\left(\D^{(k-1)/2}u_i\right)\\ \no
&=&
\int_M \left|\nabla\left( \D^{(k-1)/2}u_i\right)\right|^2
\\ \no &\geq& 0.
\en
Thus the inequality at the left hand side of (2.6) holds.
Observe that when $k$ is even, we have
\be
\int_M u_i(-\Delta )^k u_i&=&\int_M \Delta^{k/2-1}u_i \Delta\left(\Delta^{k/2}u_i\right)
\\ \no &= &  -\int_M \nabla \left(\Delta^{k/2-1}u_i\right)\nabla\left(\Delta^{k/2}u_i\right)
\\ \no
&\leq& \left(\int_M\left|\nabla \left(\Delta^{k/2-1}u_i\right)\right|^2\right)^{1/2}\left(\int_M\left|\nabla\left(\Delta^{k/2}u_i\right)\right|^2\right)^{1/2}\\ \no
&=& \left(-\int_M\Delta^{k/2-1}u_i\Delta^{k/2}u_i\right)^{1/2}\left(-\int_M\Delta^{k/2}u_i\Delta^{k/2+1}u_i\right)^{1/2}\\ \no
&=& \left(\int_M u_i(-\Delta)^{k-1}u_i\right)^{1/2}\left(\int_M u_i(-\Delta)^{k+1}u_i\right)^{1/2}.
\en
On the other hand, when $n$ is odd, it holds
\be
\int_M u_i(-\Delta )^k u_i&=&\int_M (-\Delta)^{(k-1)/2}u_i (-\Delta)^{(k+1)/2}u_i
\\ \no
&\leq& \left(\int_M \left((-\Delta)^{(k-1)/2}u_i\right)^2\right)^{1/2}\left(\int_M\left(\left(-\Delta \right)^{(k+1)/2}u_i\right)^2\right)^{1/2}\\ \no
&=& \left(\int_M u_i(-\Delta)^{k-1}u_i\right)^{1/2}\left(\int_M u_i(-\Delta)^{k+1}u_i\right)^{1/2}.
\en
Thus we always have
\be
\int_M u_i(-\Delta )^k u_i\leq  \left(\int_M u_i(-\Delta)^{k-1}u_i\right)^{1/2}\left(\int_M u_i(-\Delta)^{k+1}u_i\right)^{1/2}.
\en
When $k=1$ or $l=2$, the right hand side of (2.6) holds obviously. Now we consider the case that $l>2$ and $k\geq 2$.
 We claim now that  for any $k=2,\cdots, l-1$, it holds
 \be
 \left(\int_M u_i(-\Delta )^k u_i\right)^{k}\leq \left(\int_M u_i(-\Delta )^{k+1} u_i\right)^{k-1}.
\en
Since \be \no \int_{\om} u_i\Delta^2u_i =\int_{\om}
\Delta u_i\Delta u_i =-\int_{\om} \na \Delta u_i\na u_i,\en
 we have from
Schwarz inequality that
 \be  \left(\int_{\om} u_i\Delta^2 u_i\right)^2
 \leq
 \left(\int_{\om} |\na \D u_i|^2\right)\left(\int_{\om} |\na u_i|^2\right)
= -\int_{\om} \D u_i\D^2 u_i = \int_{\om} u_i(-\D^3 u_i). \en
Hence (2.11) holds when $k=2$. Suppose that
(2.11) holds for $k-1$, that is
 \be
\left(\int_{\om} u_i(-\Delta )^{k-1} u_i\right)^{k-1}\leq
\left(\int_{\om} u_i(-\Delta  )^{k} u_i\right)^{k-2}. \en
 Substituting (2.13) into (2.10),
we know that (2.6) is true for $k$. Using (2.6) repeatedly, we get
\be \no \int_{\om} u_i(-\Delta )^k u_i\leq \left(\int_{\om} u_i(-\Delta
)^{k+1}u_i\right)^{(k-1)/k} \leq \cdots \leq \left(\int_{\om}
u_i(-\Delta )^{l} u_i\right)^{(k-1)/(l-1)}=\laa_i^{(k-1)/(l-1)}.
\en
This completes the proof of Lemma 2.1.

\vskip0.3cm {\bf Lemma 2.2.}  Let $\{a_i\}_{i=1}^m$,
$\{b_i\}_{i=1}^m$ and $\{c_i\}_{i=1}^m$ be three  sequences of
non-negative real numbers with $\{a_i\}$ decreasing  and $\{b_i\}$
and $\{c_i\}_{i=1}^m$ increasing. Then the following inequality
holds: \be \left(\sum_{i=1}^m a_i^2b_i\right)\left(\sum_{i=1}^m a_i
c_i\right)\leq \left(\sum_{i=1}^m a_i^2 \right)\left(\sum_{i=1}^m
a_i b_i c_i\right). \en
{\it Proof.} When $m=1$, (2.14) holds
trivally. Suppose that (2.14) holds when  $m=k$, that is \be
\left(\sum_{i=1}^k a_i^2b_i\right)\left(\sum_{i=1}^k a_i
c_i\right)\leq \left(\sum_{i=1}^k a_i^2\right)\left(\sum_{i=1}^k a_i
b_i c_i\right). \en Then when $m=k+1$, we have from (2.15) that \be
& & \left(\sum_{i=1}^{k+1} a_i^2\right) \left(\sum_{i=1}^{k+1} a_i
b_ic_i\right)-\left(\sum_{i=1}^{k+1} a_i^2b_i\right)\left(\sum_{i=1}^{k+1} a_i c_i\right) \\
\no &=& \left(\sum_{i=1}^{k} a_i^2\right) \left(\sum_{i=1}^k a_i b_i
c_i\right)-\left(\sum_{i=1}^{k} a_i^2b_i\right)\left(\sum_{i=1}^k
a_i c_i\right) + a_{k+1}^2\sum_{i=1}^k a_i b_i c_i
 \\ \no & & -a_{k+1}^2 b_{k+1}\sum_{i=1}^ka_i c_i+a_{k+1}b_{k+1}c_{k+1}\sum_{i=1}^k a_i^2-a_{k+1}c_{k+1}\sum_{i=1}^k a_i^2b_i\\ \no & \geq&  a_{k+1}^2\sum_{i=1}^k a_i b_i c_i
 -a_{k+1}^2 b_{k+1}\sum_{i=1}^ka_i c_i+a_{k+1}b_{k+1}c_{k+1}\sum_{i=1}^k a_i^2-a_{k+1}c_{k+1}\sum_{i=1}^k a_i^2b_i\\ \no &=&
 -a_{k+1}^2\sum_{i=1}^k(b_{k+1}-b_i)a_ic_i+a_{k+1}c_{k+1}\sum_{i=1}^k a_i^2(b_{k+1}-b_i)\\ \no
 &=&\sum_{i=1}^k a_{k+1}a_i(b_{k+1}-b_i)(c_{k+1}a_i-a_{k+1}c_i)
 \\ \no &\geq&0.
 \en
 Where in the last inequality we have used the fact that
 \be\no
  a_{k+1}a_i(b_{k+1}-b_i)(c_{k+1}a_i-a_{k+1}c_i)\geq 0, \ \ i=1,\cdots,k.
 \en
 Thus (2.14) holds for $m=k+1$. This completes the proof of Lemma 2.2.

\vskip0.3cm
The following result is the so called {\it Reverse Chebyshev Inequality} (cf. [HLP]).
\vs
{\bf Lemma 2.4.} Suppose
$\{a_i\}_{i=1}^m$ and $\{b_i\}_{i=1}^m$ are two real sequences with
$\{a_i\}$ increasing  and $\{b_i\}$ decreasing. Then the following
inequality holds: \be \sum_{i=1}^m a_ib_i\leq \fr 1m
\left(\sum_{i=1}^m a_i\right)\left(\sum_{i=1}^m b_i\right). \en

\vs
We are now ready to prove the main results in this paper.
\vskip0.3cm
{\it Proof of Theorem 1.1.} With the notations as above, we consider now the special case that $\om$ is a connected bounded domain in
$I\!\! R^{n}$. Denote by $x_1,\cdots, x_n$ the coordinate functions on $I\!\! R^{n}$ and let us
 decompose the vector-valued functions $x_{\alpha}\nabla u_i$ as
\be
x_{\alpha}\nabla u_i=\nabla h_{\alpha i}+ W_{\alpha i},
\en
where
$h_{\alpha i}\in H_{2,D}^l(\om),$ $\nabla h_{\alpha i}$ is the projection of $x_{\alpha} \nabla u_i$ in ${\bf H}_{1, D}^2(\om)$ and $W_{\alpha i}\ \bot\ {\bf H}_{1, D}^2(\om)$. Thus we have
\be
W_{\alpha i}|_{\pa\om}\ =\ 0, \ \ {\it and} \ \ (W_{\alpha i}, \nabla u)=\int_{\om}
\langle W_{\alpha i}, \nabla u\rangle=0, \ \ {\rm for\ any} \ \ u\in H_{l,D}^2(\om)
\en
and from the discussions in [CY3] and [WX2] we know that
\be
{\rm div}\ W_{\alpha i}=0,
\en
where for a vector field $Z$ on $\om$, ${\rm div}\ Z$ denotes the divergence of $Z$.

For each $\alpha=1,\cdots, n$, $i=1,\cdots, k$, consider the  functions $\phi_{\alpha i}: \om\rightarrow I\!\!R$, given by
\be
 \phi_{\alpha i}=h_{\alpha i}-\sum_{j=1}^ka_{\alpha ij}u_j,
 \en
 where
 \be
 a_{\alpha ij}=\int_{\om}x_{\alpha}\lan\nabla u_i, \nabla u_j\ra=a_{\alpha ji}.
\en
We have
\be
\phi_{\alpha i}|_{\pa \om}=\left.\fr{\pa \phi_{\alpha i}}{\pa \nu}\right|_{\pa \om}=\cdots\left. \fr{\pa^{l-1} \phi_{\alpha i}}{\pa \nu^{l-1}}\right|_{\pa \om}=   0,
\en
\be
(\phi_{\alpha i}, u_j)_D= \int_{\om}\langle\nabla \phi_{\alpha i}, \nabla u_j\rangle=0, \ \ \forall j=1,\cdots, k.
\en
It then follows from the Rayleigh-Ritz inequality for $\laa_{k+1}$ that
\be
\laa_{k+1}\int_{\om}|\nabla \phi_{\alpha i}|^2\leq \int_D\phi_{\alpha i}(-\Delta)^l\phi_{\alpha i}, \ \ \forall\alpha =1,\cdots, n, \ \ i=1,\cdots, k.
\en
Since ${\rm div}\ W_{\alpha i}=0$, we have from (2.18) and (2.21)  that
\be\no
\Delta \phi_{\alpha i}&=&\Delta h_{\alpha i}-\sum_{j=1}^ka_{\alpha ij}\Delta u_j\\ \no
&=& {\rm div}(x_{\alpha}\nabla u_i)-\sum_{j=1}^ka_{\alpha ij}\Delta u_j
\\ \no &=&  u_{i,{\alpha}}+x_{\alpha}\Delta u_i-\sum_{j=1}^ka_{\alpha ij}\Delta u_j,
\en
where $u_{i,{\alpha}}=\fr{\pa u_i}{\pa x_{\alpha}}$.
Thus we have
\be
(-\D)^l\phi_{\alpha i}&=&(-1)^l\D^{l-1}\left(u_{i,{\alpha}}+x_{\alpha}\Delta u_i\right)+\sum_{j=1}^ka_{\alpha ij}\laa_j\Delta u_j.
\en
Since
$$
\int_{\om} \phi_{\alpha i}\Delta u_j=-\int_{\om}\langle \nabla \phi_{\alpha i}, \nabla u_j\rangle=0,
$$
\be\no
\D^{l-2}(x_{\alpha}\D u_i)=2(l-2)(\D^{l-2}u_i)_{,\alpha} +x_{\alpha}\D^{l-1} u_i,
\en
We have
\be
& &\int_{\om}\phi_{\alpha i}(-\Delta)^l\phi_{\alpha i}\\ \no
&=&\int_{\om}\phi_{\alpha i}(-1)^l\Delta^{l-1}\left(u_{i,{\alpha}}+x_{\alpha}\Delta u_i\right)
\\ \no
&=&\int_{\om}h_{\alpha i}(-1)^l\Delta^{l-1}\left(u_{i,{\alpha}}+x_{\alpha}\Delta u_i\right)
-\sum_{j=1}^ka_{\alpha ij}\int_{\om}u_j(-\Delta)^{l}h_{{\alpha}i}\\ \no
&=&\int_{\om}\D h_{\alpha i}(-1)^l\Delta^{l-2}\left(u_{i,{\alpha}}+x_{\alpha}\Delta u_i\right)
-\sum_{j=1}^ka_{\alpha ij}\int_{\om}h_{{\alpha}i}(-\Delta)^{l}u_j\\ \no
&=&\int_{\om}\D h_{\alpha i}(-1)^l((\Delta^{l-2}u_i)_{,{\alpha}}+\Delta^{l-2}(x_{\alpha}\Delta u_i))
+\sum_{j=1}^k\laa_ja_{\alpha ij}\int_{\om}h_{{\alpha}i}\D u_j\\ \no
&=&\int_{\om}(-1)^l(u_{i,{\alpha}}+x_{\alpha}\Delta u_i)((\Delta^{l-2}u_i)_{,{\alpha}}+\Delta^{l-2}(x_{\alpha}\Delta u_i))
-\sum_{j=1}^k\laa_ja_{\alpha ij}\int_{\om} \lan \na h_{\alpha i}, \na u_j\ra
\\ \no
&=&\int_{\om}(-1)^l(u_{i,{\alpha}}+x_{\alpha}\Delta u_i)((2l-3)(\Delta^{l-2}u_i)_{,{\alpha}}+x_{\alpha}\Delta^{l-1} u_i))
-\sum_{j=1}^k\laa_ja_{\alpha ij}\int_{\om} \lan \na h_{\alpha i}, \na u_j\ra
 \\ \no
&=&\int_{\om}(-1)^l((2l-3)(u_{i,{\alpha}}(\Delta^{l-2}u_i)_{,{\alpha}}+x_{\alpha}\Delta u_i(\Delta^{l-2}u_i)_{,{\alpha}})+ u_{i,{\alpha}}x_{\alpha}\Delta^{l-1}u_i+x_{\alpha}^2\D u_i\Delta^{l-1} u_i)
-\sum_{j=1}^k\laa_ja_{\alpha ij}^2
\en
Let us make some calculations.
Since
\be\no
\D^{l-1}(x_{\alpha}u_i)=2(l-1)(\D^{l-2}u_i)_{,\alpha} +x_{\alpha}\D^{l-1} u_i,
\en
we have
\be
\int_{\om}x_{\alpha}u_i (\D^{l-1}u_i)_{,\alpha}  &=&\int_{\om}x_{\alpha}u_i \D^{l-1} u_{i,\alpha}\\ \no &=&\int_{\om}\D^{l-1}(x_{\alpha}u_i) u_{i,\alpha}\\ \no &=&\int_{\om}\left(2(l-1)(\D^{l-2}u_i)_{,\alpha} +x_{\alpha}\D^{l-1} u_i\right) u_{i,\alpha}.
\en
On the other hand, it holds
\be
\int_{\om}x_{\alpha}u_i (\D^{l-1}u_i)_{,\alpha} =-\int_{\om}\D^{l-1}u_i(u_i+x_{\alpha} u_{i,\alpha}).
\en
Combining  (2.28) and (2.29), we obtain
\be & &
\int_{\om}x_{\alpha}u_i(\D^{l-1}u_i)_{,\alpha}
 \\ \no &=&
 \int_M\left\{(l-1)(\D^{l-2}u_i)_{,\alpha}u_{i,\alpha}-\fr 12u_i\D^{l-1}u_i \right\}
\en
Hence
\be
\int_{\om}x_{\alpha}u_{i,\alpha}\D^{l-1} u_i&=&-\int_{\om}u_i(\D^{l-1} u_i+x_{\alpha}(\D^{l-1}u_i)_{,\alpha})\\ \no
&=& -\int_M\left\{(l-1)(\D^{l-2}u_i)_{,\alpha}u_{i,\alpha}+\fr 12u_i\D^{l-1}u_i \right\}
\en
and consequently, we have
\be
\int_{\om}x_{\alpha}\Delta u_i(\Delta^{l-2}u_i)_{,{\alpha}}
&=&\int_{\om}x_{\alpha}\Delta u_i\Delta^{l-2}u_{i,{\alpha}}\\ \no
&=&\int_{\om}\Delta^{l-2}(x_{\alpha}\Delta u_i)u_{i,{\alpha}}\\ \no
&=&\int_{\om}u_{i,\alpha}\left(2(l-2)(\D^{l-2}u_i)_{,\alpha} +x_{\alpha}\D^{l-1} u_i\right)
\\ \no &=&\int_{\om}\left\{(l-3)(\D^{l-2}u_i)_{,\alpha}u_{i,\alpha}-\fr 12u_i\D^{l-1}u_i\right\}.
\en
Also, one has
\be
\int_{\om}u_ix_{\alpha}^2\D u_i&=&-\int_{\om}x_{\alpha}^2|\na u_i|^2-2\int_{\om}x_{\alpha}u_i u_{i,\alpha}\\
\no
&=&-\int_{\om}x_{\alpha}^2|\na u_i|^2+ \int_{\om}u_i^2,
\en
\be
\int_{\om} x_{\alpha}^2\D u_i\Delta^{l-1} u_i&=&\int_{\om}u_i\D(x_{\alpha}^2\Delta^{l-1} u_i)
\\ \no
&=&\int_{\om}u_i(2\D^{l-1}u_i+x_{\alpha}^2\D^lu_i+4x_{\alpha}(\D^{l-1}u_i)_{,\alpha})
\\ \no
&=&\int_{\om}u_i(2\D^{l-1}u_i+(-1)^{l-1}\laa_i x_{\alpha}^2\D u_i+4x_{\alpha}(\D^{l-1}u_i)_{,\alpha}).
\en
Combining (2.30), (2.33) and (2.34), we get
\be& &
\int_{\om} x_{\alpha}^2\D u_i\Delta^{l-1} u_i\\ \no &=&4(l-1)\int_{\om}(\D^{l-2}u_i)_{,\alpha}u_{i,\alpha}
+(-1)^{l-1}\laa_i\left\{-\int_{\om}x_{\alpha}^2|\na u_i|^2+\int_{\om}u_i^2\right\}.
\en
Substituting (2.32), (2.33) and (2.35) into (2.27), one gets
\be
\int_{\om}\phi_{\alpha i}(-\Delta)^l\phi_{\alpha i}
&=&\int_{\om}(-1)^l\left\{(-l+1)u_i\D^{l-1}u_i+(2l^2-4l+3)(\D^{l-2}u_i)_{,\alpha}u_{i,\alpha}\right\}
\\ \no & & +\laa_i\left\{\int_{\om}x_{\alpha}^2|\na u_i|^2- \int_{\om}u_i^2\right\}-\sum_{j=1}^k\laa_ja_{\alpha ij}^2.
\en
It is easy to see that
\be
||x_{\alpha}\nabla u_i||^2=||\nabla h_{\alpha i}||^2+||W_{\alpha i}||^2, \ \ ||\nabla h_{\alpha i}||^2=||\nabla \phi_{\alpha i}||^2+\sum_{j=1}^ka_{\alpha ij}^2,
\en
where for a vector field $Z$ on $\om$, $||Z||^2=\int_{\om} |Z|^2$.
Combining (2.25), (2.36) and (2.37), we infer
\be
(\laa_{k+1}-\laa_i)
||\nabla \phi_{\alpha i}||^2&\leq& \int_{\om}(-1)^l\left\{(-l+1)u_i\D^{l-1}u_i+(2l^2-4l+3)(\D^{l-2}u_i)_{,\alpha}u_{i,\alpha}\right\}
\\ \no & & -\laa_i(||u_i||^2-||W_{\alpha i}||^2)+\sum_{j=1}^k(\laa_i-\laa_j)a_{\alpha ij}^2,
\en
Observe that $\na(x_{\alpha}u_i)=u_i\na x_{\alpha}+ x_{\alpha}\na u_i\in {\bf H}_{1, D}^2(\om )$. For $A_{\alpha i}=\na(x_{\alpha} u_i-h_{\alpha i})$, we have
\be
u_i\na x_{\alpha}=A_{\alpha i}-W_{\alpha i}
\en
and so
$$||u_{i}||^2=||u_i \na x_{\alpha}||^2=||W_{\alpha i}||^2+||A_{\alpha i}||^2.$$
Because of $(\na u_{i,\alpha}, W_{\alpha i})=0$, it follows that
\be\no
2||u_{i,\alpha}||^2=-2\int_{\om}\lan u_i\na x_{\alpha}, \na u_{i,\alpha}\ra=-2\int_{\om}\lan A_{\alpha i}, \na u_{i,\alpha}\ra\leq \laa_i^{1/(l-1)}||A_{\alpha i}||^2+\fr 1{\laa_i^{1/(l-1)}}|| \na u_{i,\alpha}||^2
\en
which gives
\be
-\laa_i||A_{\alpha i}||^2\leq -2\laa_i^{(l-2)/(l-1)}||u_{i,\alpha}||^2 +\laa_i^{(l-3)/(l-1)}|| \na u_{i,\alpha}||^2
\en
Introducing (2.40) into (2.38), we get
\be
(\laa_{k+1}-\laa_i)
||\nabla \phi_{\alpha i}||^2&\leq& \int_{\om}(-1)^l\left\{(-l+1)u_i\D^{l-1}u_i+(2l^2-4l+3)(\D^{l-2}u_i)_{,\alpha}u_{i,\alpha}\right\}
\\ \no & & -2\laa_i^{(l-2)/(l-1)}||u_{i,\alpha}||^2 +\laa_i^{(l-3)/(l-1)}|| \na u_{i,\alpha}||^2 +\sum_{j=1}^k(\laa_i-\laa_j)a_{\alpha ij}^2,
\en
Since
\be\no
-2\int_{\om} x_{\alpha}\lan \na u_i, \na u_{i,\alpha}\ra &=&2\int_{\om}u_{i,\alpha}^2 + 2\int_{\om}x_{\alpha}u_{i,\alpha}\D u_i\\ \no
&=& 2\int_{\om}u_{i,\alpha}^2 + 2\int_{\om} u_i \D(x_{\alpha}u_{i,\alpha})\\ \no
&=& 2\int_{\om}u_{i,\alpha}^2 + 2\int_{\om} u_i x_{\alpha}(\D u_i)_{, \alpha}+4\int_{\om} \lan u_i \na x_{\alpha} , \na u_{i,\alpha}\ra
\\ \no
&=& 2\int_{\om}u_{i,\alpha}^2 - 2\int_{\om}\D u_i(u_i+ x_{\alpha} u_{i,\alpha})-4\int_{\om} u_{i,\alpha}{\rm div}(u_i \na x_{\alpha})\\ \no
&=& 2\int_{\om}u_{i,\alpha}^2 +2- 2\int_{\om}x_{\alpha} u_{i,\alpha}\D u_i-4\int_{\om}u_{i,\alpha}^2\\ \no
&=& -2\int_{\om}u_{i,\alpha}^2 +2+2\int_{\om}\lan \na u_i, \na(x_{\alpha} u_{i,\alpha})\ra \\ \no
&=& 2+2\int_{\om} x_{\alpha}\lan \na u_i, \na u_{i,\alpha}\ra,
\en
we have
\be
-2\int_{\om} x_{\alpha}\lan \na u_i, \na u_{i,\alpha}\ra=1.
\en
Set
\be\no
d_{\alpha ij}=\int_{\om}\lan \na u_{i,\alpha}, \na u_j\ra ;
\en
then $d_{\alpha ij}=-d_{\alpha ji}$ and we have from (2.18), (2.20) and (2.21) that
\be
1&=&-2\int_{\om} x_{\alpha}\lan \na u_i, \na u_{i,\alpha}\ra\\ \no
&=& -2\int_{\om} \lan \na h_{\alpha i}, \na u_{i,\alpha}\ra
\\ \no &=& -2\int_{\om} \lan \na \phi_{\alpha i}, \na u_{i,\alpha}\ra -2\sum_{j=1}^k a_{\alpha ij}d_{\alpha ij}.
\en
Thus, we have
\be
& & (\laa_{k+1}-\laa_i)^2\left(1+2\sum_{j=1}^k a_{\alpha ij}d_{\alpha ij}\right)\\ \no
&=&(\laa_{k+1}-\laa_i)^2\left(-2\na \phi_{\alpha i}, \left(\na u_{i,\alpha}-\sum_{j=1}^k d_{\alpha ij}\na u_j\right)\right)\\ \no &\leq&\delta (\laa_{k+1}-\laa_i)^3||\na \phi_{\alpha i}||^2+\fr 1\delta (\laa_{k+1}-\laa_i)\left(||\na u_{i,\alpha}||^2-\sum_{j=1}^k d_{\alpha ij}^2\right),
\en
where $\delta$ is any positive constant.
Substituting (2.41) into (2.44), we get
\be\no
& & (\laa_{k+1}-\laa_i)^2\left(1+2\sum_{j=1}^k a_{\alpha ij}d_{\alpha ij}\right)\\ \no
&\leq&\delta (\laa_{k+1}-\laa_i)^2\left(\int_{\om}(-1)^l\left\{(-l+1)u_i\D^{l-1}u_i+(2l^2-4l+3)(\D^{l-2}u_i)_{,\alpha}u_{i,\alpha}\right\}
\right.
\\ \no & & \left. -2\laa_i^{(l-2)/(l-1)}||u_{i,\alpha}||^2 +\laa_i^{(l-3)/(l-1)}|| \na u_{i,\alpha}||^2 +\sum_{j=1}^k(\laa_i-\laa_j)a_{\alpha ij}^2\right)
\\ \no & & +\fr 1\delta (\laa_{k+1}-\laa_i)\left(||\na u_{i,\alpha}||^2-\sum_{j=1}^k d_{\alpha ij}^2\right),
\en
Summing on $i$ from $1$ to $k$  and noticing the fact that $a_{\alpha ij}=a_{\alpha ji}, d_{\alpha ij}=-d_{\alpha ji}$, we infer
\be
& & \sum_{i=1}^k(\laa_{k+1}-\laa_i)^2-2\sum_{i, j=1}^k (\laa_{k+1}-\laa_i)(\laa_i-\laa_j)a_{\alpha ij}d_{\alpha ij}\\ \no
&\leq&\delta\left(\sum_{i=1}^k (\laa_{k+1}-\laa_i)^2\left(\int_{\om}(-1)^l\left\{(-l+1)u_i\D^{l-1}u_i+(2l^2-4l+3)(\D^{l-2}u_i)_{,\alpha}u_{i,\alpha}\right\}
\right.\right.
\\ \no & & \left. \left.-2\laa_i^{(l-2)/(l-1)}||u_{i,\alpha}||^2 +\laa_i^{(l-3)/(l-1)}|| \na u_{i,\alpha}||^2\right) -\sum_{i, j=1}^k(\laa_{k+1}-\laa_i)(\laa_i-\laa_j)^2a_{\alpha ij}^2\right)\\ \no & &
+\fr 1\delta\left(\sum_{i=1}^k  (\laa_{k+1}-\laa_i)||\na u_{i,\alpha}||^2-\sum_{i, j=1}^k  (\laa_{k+1}-\laa_i)d_{\alpha ij}^2\right),
\en
which gives
\be
& & \sum_{i=1}^k(\laa_{k+1}-\laa_i)^2\\ \no
&\leq&\delta\sum_{i=1}^k (\laa_{k+1}-\laa_i)^2\left(\int_{\om}(-1)^l\left\{(-l+1)u_i\D^{l-1}u_i+(2l^2-4l+3)(\D^{l-2}u_i)_{,\alpha}u_{i,\alpha}\right\}
\right.
\\ \no & & \left.-2\laa_i^{(l-2)/(l-1)}||u_{i,\alpha}||^2 +\laa_i^{(l-3)/(l-1)}|| \na u_{i,\alpha}||^2\right)
+\fr 1\delta \sum_{i=1}^k  (\laa_{k+1}-\laa_i)||\na u_{i,\alpha}||^2,
\en
Taking sum for $\alpha $ from $1$ to $n$, we get
\be\no
& & n\sum_{i=1}^k(\laa_{k+1}-\laa_i)^2\\ \no
&\leq&\delta\sum_{i=1}^k (\laa_{k+1}-\laa_i)^2\left(\int_{\om}(-1)^l\left\{n(-l+1)u_i\D^{l-1}u_i+(2l^2-4l+3)\lan\na(\D^{l-2}u_i),
\na u_i\ra \right\}
\right.
\\ \no & & \left.-2\laa_i^{(l-2)/(l-1)} +\laa_i^{(l-3)/(l-1)}\sum_{\alpha=1}^n || \na u_{i,\alpha}||^2\right)
+\fr 1\delta \sum_{i=1}^k  (\laa_{k+1}-\laa_i)\sum_{\alpha =1}^n||\na u_{i,\alpha}||^2
\\ \no
&=&\delta\sum_{i=1}^k (\laa_{k+1}-\laa_i)^2\left( -2\laa_i^{(l-2)/(l-1)} +\laa_i^{(l-3)/(l-1)}\sum_{\alpha=1}^n || \na u_{i,\alpha}||^2 +(2l^2+(n-4)l+3-n) \int_{\om}u_i(-\D)^{l-1}u_i
\right)
\\ \no & &
+\fr 1\delta \sum_{i=1}^k  (\laa_{k+1}-\laa_i)\sum_{\alpha =1}^n||\na u_{i,\alpha}||^2.
\en
But
\be\no
\sum_{\alpha =1}^k||\na u_{i,\alpha}||^2&=&-\int_{\om}\sum_{\alpha =1}^k u_{i,\alpha}\D  u_{i,\alpha}\\ \no
&=&-\int_{\om}\sum_{\alpha =1}^k u_{i,\alpha}(\D  u_i)_{,\alpha}
\\ \no &=&
\int_{\om}\sum_{\alpha =1}^k u_{i,\alpha\alpha }\D  u_i
\\ \no &=&
\int_{\om}(\D  u_{i})^2 \\ \no &=&
\int_{\om}u_i\D^2  u_{i},
\en
where $u_{i,\alpha\alpha }=\fr{\pa^2 u_i}{\pa x_{\alpha}^2}$.
Thus, we have
\be
& & n\sum_{i=1}^k(\laa_{k+1}-\laa_i)^2
\\ \no
&\leq&\delta\sum_{i=1}^k (\laa_{k+1}-\laa_i)^2\left( -2\laa_i^{(l-2)/(l-1)} +\laa_i^{(l-3)/(l-1)}\int_{\om}u_i\D^2  u_i\right.\\ \no
& &\left. +(2l^2+(n-4)l+3-n) \int_{\om}u_i(-\D)^{l-1}u_i
\right)
+\fr 1\delta \sum_{i=1}^k  (\laa_{k+1}-\laa_i)\int_{\om}u_i\D^2  u_{i}.
\en
 Taking $k=2$ and $k=l-1$ in (2.6), respectively, one gets
\be\no
\int_{\om}u_i(-\D)^{l-1}u_i\leq \laa_i^{(l-2)/(l-1)}, \ \ \ \int_{\om}u_i\D^2u_i\leq \laa_i^{1/(l-1)}
\en
which, combining with (2.47) implies that
\be
& &\no n\sum_{i=1}^k(\laa_{k+1}-\laa_i)^2\\ \no
&\leq&\delta (2l^2+(n-4)l+2-n) \sum_{i=1}^k  (\laa_{k+1}-\laa_i)^2 \laa_i^{(l-2)/(l-1)}
+\fr 1\delta \sum_{i=1}^k  (\laa_{k+1}-\laa_i)\laa_i^{1/(l-1)}.
\en
Taking
\be
\no
\delta=\fr{\left\{\sum_{i=1}^k  (\laa_{k+1}-\laa_i)\laa_i^{1/(l-1)}\right\}^{1/2}}{\left\{(2l^2+(n-4)l+2-n)\sum_{i=1}^k  (\laa_{k+1}-\laa_i)^2 \laa_i^{(l-2)/(l-1)}\right\}^{1/2}},
\en
we get (1.13). This completes the proof of Theorem 1.1.
\vskip0.3cm
{\it Proof of Corollary 1.1.}  It follows from (2.17) that
\be
 \sum_{i=1}^k  (\laa_{k+1}-\laa_i)\laa_i^{1/(l-1)}\leq \fr 1k\left(\sum_{i=1}^k  (\laa_{k+1}-\laa_i)\right)\left(\sum_{i=1}^k
\laa_i^{1/(l-1)}\right)
\en
and
\be
\sum_{i=1}^k  (\laa_{k+1}-\laa_i)^2 \laa_i^{(l-2)/(l-1)}\leq\fr 1k\left(\sum_{i=1}^k  (\laa_{k+1}-\laa_i)^2\right)\left(\sum_{i=1}^k\laa_i^{(l-2)/(l-1)}\right)
\en
Introducing (2.48) and (2.49) into (1.13), we infer
\be\no
& &\sum_{i=1}^k(\laa_{k+1}-\laa_i)^2\\ \no
&\leq& \fr{ 4(2l^2+(n-4)l+2-n)}{k^2n^2}\left( \sum_{i=1}^k  (\laa_{k+1}-\laa_i)\right)\left(\sum_{i=1}^k
\laa_i^{1/(l-1)}\right)\left(\sum_{i=1}^k\laa_i^{(l-2)/(l-1)}\right).
\en
Solving this quadratic polynomial about $\laa_{k+1}$, one gets (1.14).

From (2.14), we have
\be\no
\left\{ \sum_{i=1}^k  (\laa_{k+1}-\laa_i)^2 \laa_i^{(l-2)/(l-1)}\right\}
\left\{ \sum_{i=1}^k  (\laa_{k+1}-\laa_i)\laa_i^{1/(l-1)}\right\}
\leq \left\{ \sum_{i=1}^k  (\laa_{k+1}-\laa_i)^2\right\}\left\{ \sum_{i=1}^k  (\laa_{k+1}-\laa_i)\laa_i\right\}.
\en
It then follows from (1.13) that
\be\no
\sum_{i=1}^k  (\laa_{k+1}-\laa_i)^2\leq \fr{ 4(2l^2+(n-4)l+2-n)}{n^2}\sum_{i=1}^k  (\laa_{k+1}-\laa_i)\laa_i,
\en
which implies (1.15). This completes the proof of Corollary 1.1.

\vskip0.3cm
{\it Proof of Theorem 1.2.} We use  the same notations as in the beginning of this section and take $M$ to be the unit $n$-sphere $S^n(1)$.
Let $x_1, x_2,\cdots, x_{n+1}$ be the standard  coordinate functons of the Euclidean space  $I\!\!R^{n+1}$; then
$$S^n(1)=\left\{(x_1,\dots, x_{n+1})\in  I\!\!R^{n+1}; \sum_{\alpha=1}^{n+1}x_{\alpha}^2=1\right\}.$$ It is well known that
\be
\Delta x_{\alpha}=-nx_{\alpha},\ \ \  \alpha =1,\cdots, n+1.
\en
As in the proof of Theorem 1.1, we
 decompose the vector-valued functions $x_{\alpha}\nabla u_i$ as
\be
x_{\alpha}\nabla u_i=\nabla h_{\alpha i}+ W_{\alpha i},
\en
where
$h_{\alpha i}\in H_{l,D}^2(\om),$ $\nabla h_{\alpha i}$ is the projection of $x_{\alpha} \nabla u_i$ in ${\bf H}_{1, D}^2(\om)$, $W_{\alpha i}\ \bot\ {\bf H}_{1, D}^2(\om)$
 and
\be
W_{\alpha i}|_{\pa\om}\ =\ 0, \ \  {\rm div\ } W_{\alpha i}=0.
\en
We also consider the  functions $\phi_{\alpha i}: \om\rightarrow I\!\!R$, given by
\be
 \phi_{\alpha i}=h_{\alpha i}-\sum_{j=1}^kb_{\alpha ij}u_j, \ \
 b_{\alpha ij}=\int_{\om}x_{\alpha}\lan\nabla u_i, \nabla u_j\ra=b_{\alpha ji}.
\en
 Then
\be\no
\phi_{\alpha i}|_{\pa \om}=\left.\fr{\pa \phi_{\alpha i}}{\pa \nu}\right|_{\pa \om}=\cdots\left. \fr{\pa^{l-1} \phi_{\alpha i}}{\pa \nu^{l-1}}\right|_{\pa \om}=   0,
\en
\be\no
(\phi_{\alpha i}, u_j)_D= \int_{\om}\langle\nabla \phi_{\alpha i}, \nabla u_j\rangle=0, \ \ \forall j=1,\cdots, k
\en
and we have the basic Rayleigh-Ritz inequality for $\laa_{k+1}$ :
\be
\laa_{k+1}\int_{\om}|\nabla \phi_{\alpha i}|^2\leq \int_D\phi_{\alpha i}(-\Delta)^l\phi_{\alpha i}, \ \ \forall\alpha =1,\cdots, n, \ \ i=1,\cdots, k.
\en
We have
\be
\Delta \phi_{\alpha i}= \langle \nabla x_{\alpha}, \nabla u_i\rangle+x_{\alpha}\Delta u_i-\sum_{j=1}^kb_{\alpha ij}\Delta u_j
\en
and as in the proof of (2.27),
\be
& & \int_{\om}\phi_{\alpha i}(-\Delta)^l\phi_{\alpha i}\\ \no & &
=\int_{\om}(-1)^l(\langle \nabla x_{\alpha}, \nabla u_i\rangle +x_{\alpha}\Delta u_i)\Delta^{l-2}(\langle \nabla x_{\alpha}, \nabla u_i\rangle+x_{\alpha}\Delta u_i)
-\sum_{j=1}^k\laa_jb_{\alpha ij}^2.
\en
For a function $g$ on $\om$, we have (cf. (2.31) in [WX2])
\be
\Delta \langle \nabla x_{\alpha}, \nabla g\rangle = -2x_{\alpha}\Delta g+\langle \nabla x_{\alpha}, \nabla((\Delta +n-2)g)\rangle.
\en
For each $q=0, 1,\cdots$, thanks to (2.50) and (2.57), there are polynomials
$F_q$ and $G_q$ of degree $q$ such that
\be
\Delta^{q}(\langle \nabla x_{\alpha}, \nabla u_i\rangle+x_{\alpha}\Delta u_i)
=x_{\alpha}F_q(\D)\D u_i+\lan \na x_{\alpha}, \na (G_q(\D)u_i)\ra.
\en
It is obvious that
\be
F_0=1, \ \ G_0=1.
\en
It follows from (2.50) and (2.57) that
\be
\D (x_{\alpha}\Delta u_i+\langle \nabla x_{\alpha}, \nabla u_i\rangle)=x_{\alpha}(\D-(n+2))\D u_i+\langle \nabla x_{\alpha}, \na((3\D+n-2)u_i)\ra
\en
which gives
\be
F_1(t)=t-(n+2), \ \ G_1(t)=3t+n-2.
\en
Also, when $q\geq 2$, we have
\be
& & \Delta^{q}(\langle \nabla x_{\alpha}, \nabla u_i\rangle+x_{\alpha}\Delta u_i)
\\ \no
&=& \D (\Delta^{q-1}(\langle \nabla x_{\alpha}, \nabla u_i\rangle+x_{\alpha}\Delta u_i))\\ \no &=&
\D(x_{\alpha}F_{q-1}(\D)\D u_i+\lan \na x_{\alpha}, \na (G_{q-1}(\D)u_i)\ra)
\\ \no &=& x_{\alpha}((\D-n)F_{q-1}(\D)-2 G_{q-1}(\D))\D u_i+\lan \na x_{\alpha}, \na(((\D +n-2)G_{q-1}(\D)+2\D F_{q-1}(\D ) )u_i)\ra
\en
which, combining with (2.58), implies that
\be
F_q(\D)=(\D-n)F_{q-1}(\D)-2 G_{q-1}(\D), \ \ \ q=2,\cdots,
\en
\be
G_q(\D)=(\D +n-2)G_{q-1}(\D)+2\D F_{q-1}(\D ),\ \ \ q=2,\cdots.
\en
It then follows from (2.63) and (2.64) that
\be\no
F_q(\D)&=&(\D-n)F_{q-1}(\D)-2((\D +n-2)G_{q-2}(\D)+2\D F_{q-2}(\D ))\\ \no
&=&(\D-n)F_{q-1}(\D)+(\D +n-2)(F_{q-1}(\D)-(\D -n)F_{q-2}(\D ))-4 \D F_{q-2}(\D )\\ \no &=&
(2\D-2)F_{q-1}(\D)-(\D^2+2\D -n(n-2))F_{q-2}(\D)
\en
and
\be\no
G_q(\D)&=&(\D +n-2)G_{q-1}(\D)+2\D ((\D-n)F_{q-2}(\D)-2 G_{q-2}(\D))\\ \no &=&
(\D +n-2)G_{q-1}(\D)+(\D-n)(G_{q-1}(\D)-(\D +n-2)G_{q-2}(\D))-4\D G_{q-2}(\D)\\ \no &=&
(2\D-2)G_{q-1}(\D)-(\D^2+2\D -n(n-2))G_{q-2}(\D).
\en
Thus, we have
\be& &
F_q(t)=(2t-2)F_{q-1}(t)-(t^2+2t-n(n-2))F_{q-2}(t), \\  & & G_q(t)=(2t-2)G_{q-1}(t)-(t^2+2t-n(n-2))G_{q-2}(t), \ \ q=2,\cdots.
\en
That is, the polynomials $F_q$ and $G_q$ are defined inductively by (1.17)-(1.19).
Substituting
 \be
\Delta^{l-2}(\langle \nabla x_{\alpha}, \nabla u_i\rangle+x_{\alpha}\Delta u_i)
=x_{\alpha}F_{l-2}(\D)\D u_i+\lan \na x_{\alpha}, \na (G_{l-2}(\D)u_i)\ra
\en
into (2.56), we get
\be
& &\int_{\om}\phi_{\alpha i}(-\Delta)^l\phi_{\alpha i}\\ \no &=&
\int_{\om}(-1)^l(\langle \nabla x_{\alpha}, \nabla u_i\rangle \lan \na x_{\alpha}, \na (G_{l-2}(\D)u_i)\ra
+\lan x_{\alpha}\na x_{\alpha}, \Delta u_i\na (G_{l-2}(\D)u_i)+(F_{l-2}(\D)\D u_i) \na u_i\ra)\\ \no & & +\int_{\om}(-1)^l x_{\alpha}^2\Delta u_iF_{l-2}(\D)(\D u_i) -\sum_{j=1}^k\laa_jb_{\alpha ij}^2
\en
Summing over $\alpha $  and noticing
\be
 \sum_{\alpha=1}^{n+1}x_{\alpha}^2=1, \ \  \sum_{\alpha=1}^{n+1}\langle \nabla x_{\alpha}, \nabla u_i\rangle \lan \na x_{\alpha}, \na (G_{l-2}(\D)u_i)\ra = \langle  \nabla u_i , \na (G_{l-2}(\D)u_i)\ra,
\en
 we get from (1.20) and (2.6) that
\be
& &\sum_{\alpha=1}^{n+1}\int_{\om}\phi_{\alpha i}(-\Delta)^l\phi_{\alpha i}\\ \no &=&
\int_{\om}(-1)^l\langle \nabla u_i, \na (G_{l-2}(\D)u_i)\ra
 +\int_{\om}(-1)^l \Delta u_iF_{l-2}(\D)(\D u_i) -\sum_{\alpha=1}^{n+1}\sum_{j=1}^k\laa_jb_{\alpha ij}^2
\\ \no &=&
\int_{\om}(-1)^{l-1} u_i \D (G_{l-2}(\D)u_i)
 +\int_{\om}(-1)^l  u_i\D (F_{l-2}(\D)(\D u_i)) -\sum_{\alpha=1}^{n+1}\sum_{j=1}^k\laa_jb_{\alpha ij}^2
\\ \no &=&
\int_{\om}(-1)^{l} u_i \left( \D (F_{l-2}(\D)-(G_{l-2}(\D)\right)(\D u_i)
-\sum_{\alpha=1}^{n+1}\sum_{j=1}^k\laa_jb_{\alpha ij}^2
\\ \no &=&
\int_{\om}(-1)^{l} u_i \left( \D^{l-1}+a_{l-2}\D^{l-2}+\cdots +a_1 \D +a_0\right)(\D u_i)
-\sum_{\alpha=1}^{n+1}\sum_{j=1}^k\laa_jb_{\alpha ij}^2
\\ \no &=&
\laa_i+ \int_{\om}(-1)^{l} u_i \left( a_{l-2}\D^{l-2}+\cdots +a_1 \D +a_0\right)(\D u_i)
-\sum_{\alpha=1}^{n+1}\sum_{j=1}^k\laa_jb_{\alpha ij}^2
\\ \no &\leq &
\laa_i+ \sum_{j=0}^{l-2}|a_j|\int_{\om} u_i (-\D)^{j+1}u_i-\sum_{\alpha=1}^{n+1}\sum_{j=1}^k\laa_jb_{\alpha ij}^2
\\ \no &\leq & \laa_i+ \sum_{j=0}^{l-2}|a_j|\laa_i^{ j/(l-1)}-\sum_{\alpha=1}^{n+1}\sum_{j=1}^k\laa_jb_{\alpha ij}^2.
\en
Observe from (2.51) and (2.53) that
\be
||x_{\alpha}\na u_i||^2=||\na h_{\alpha i}||^2+||W_{\alpha i}||^2=
||\na \phi_{\alpha i}||^2+||W_{\alpha i}||^2+\sum_{j=1}^kb_{\alpha ij}^2.
\en
Summing over $\alpha$, one gets
\be
1=\sum_{\alpha=1}^{n+1}\left(||\na \phi_{\alpha i}||^2+||W_{\alpha i}||^2+\sum_{j=1}^kb_{\alpha ij}^2\right).
\en
Set
\be
Z_{\alpha i}=\nabla\langle \nabla x_{\alpha},\ \nabla u_i\rangle-\fr{n-2}2x_{\alpha}\nabla u_i,\ \
c_{\alpha ij}=\int_{\om}\langle\nabla u_j,\ Z_{\alpha i}\rangle
;
\en
then
$c_{\alpha ij}=-c_{\alpha ji}$ (cf. Lemma in [WX2]). By using the same arguments as in the proof of (2.37) in [WX2], we have
\be & &
(\laa_{k+1}-\laa_i)^2\left(2||\langle\nabla x_{\alpha}, \nabla u_i\rangle||^2 +\int_{\om} \left\langle\nabla x_{\alpha}^2,\ \Delta u_i\nabla u_i\right\rangle +(n-2)|| x_{\alpha}\nabla u_i||^2+2\sum_{j=1}^k b_{\alpha ij}c_{\alpha ij}\right)
\\ \no  &\leq &\delta (\laa_{k+1}-\laa_i)^3||\nabla \phi_{\alpha i}||^2
+\fr{\laa_{k+1}-\laa_i}{\delta}\left(||Z_{\alpha i}||^2-\sum_{j=1}^kc_{\alpha ij}^2\right)
+(n-2)(\laa_{k+1}-\laa_i)^2 ||W_{\alpha i}||^2
\en
where $\delta $ is any  positive constant.
Since
\be
\sum_{\alpha=1}^{n+1}||\langle\nabla x_{\alpha}, \nabla u_i\rangle||^2=\int_{\om} |\na u_i|^2=1,
\en
we have by
summing over $\alpha$ in (2.74) from 1 to $n+1$  that
\be & &
(\laa_{k+1}-\laa_i)^2\left(n+2\sum_{\alpha=1}^{n+1}\sum_{j=1}^k b_{\alpha ij}c_{\alpha ij}\right)
\\ \no  &\leq &\delta \sum_{\alpha=1}^{n+1}(\laa_{k+1}-\laa_i)^3||\nabla \phi_{\alpha i}||^2
+\fr{\laa_{k+1}-\laa_i}{\delta}\sum_{\alpha=1}^{n+1}\left(||Z_{\alpha i}||^2-\sum_{j=1}^kc_{\alpha ij}^2\right)
\\ \no & & +(n-2)\sum_{\alpha=1}^{n+1}(\laa_{k+1}-\laa_i)^2 ||W_{\alpha i}||^2.
\en
It follows from (2.6) and (2.57) that
\be
& &
\sum_{\alpha=1}^{n+1}||\nabla\langle
\nabla x_{\alpha}, \nabla u_i\rangle||^2\\ \no &=& -\sum_{\alpha=1}^{n+1}\int_{\om}\langle
\nabla x_{\alpha}, \nabla u_i\rangle\Delta \langle
\nabla x_{\alpha}, \nabla u_i\rangle\\ \no
&=&-\sum_{\alpha=1}^{n+1}\int_{\om}\langle
\nabla x_{\alpha}, \nabla u_i\rangle\left(-2x_{\alpha}\Delta u_i+\langle \nabla x_{\alpha}, \nabla(\Delta u_i)\rangle+(n-2)\langle\nabla x_{\alpha}, \ \nabla u_i\rangle\right)
\\ \no &=& -\int_{\om}\langle\nabla u_i, \nabla(\Delta  u_i)\rangle-(n-2)||\nabla u_i||^2
\\ \no &=& \int_{\om}u_i\Delta^2  u_i-(n-2)
\\ \no &\leq& \laa_i^{1/(l-1)}-(n-2)
\en
and so
\be & &
\sum_{\alpha=1}^{n+1}||Z_{\alpha i}||^2
\\ \no &=&\int_{\om}\left|\nabla\langle
\nabla x_{\alpha}, \nabla u_i\rangle-\fr{n-2}2x_{\alpha}\nabla u_i\right|^2
\\ \no &=&\sum_{\alpha=1}^{n+1}\left(||\nabla\langle
\nabla x_{\alpha}, \nabla u_i\rangle||^2-(n-2)\int_{\om}\langle\nabla\langle
\nabla x_{\alpha}, \nabla u_i\rangle, \ x_{\alpha} \nabla u_i\rangle +\fr{(n-2)^2}4||x_{\alpha} \nabla u_i||^2\right)\\ \no &\leq&
\laa_i^{1/(l-1)}-(n-2)+(n-2)+\fr{(n-2)^2}4=\laa_i^{1/(l-1)}+\fr{(n-2)^2}4.
\en
Since
\be\no
\int_{\om}\langle
\nabla x_{\alpha}, \nabla u_i\rangle^2&=&\int_{\om}\left\langle\langle
\nabla x_{\alpha}, \nabla u_i\rangle\nabla u_i, \nabla x_{\alpha}\right\rangle
\\ \no
&=&-\int_{\om} x_{\alpha}\ {\rm div}(\langle
\nabla x_{\alpha}, \nabla u_i\rangle\nabla u_i)\\ \no
&=&-\int_{\om} \langle x_{\alpha}\nabla u_i, \nabla\langle
\nabla x_{\alpha}, \nabla u_i\rangle\rangle - \int_{\om} x_{\alpha}\langle
\nabla x_{\alpha}, \nabla u_i\rangle\Delta u_i
\\ \no
&=&-\int_{\om} \langle\nabla h_{\alpha i}+W_{\alpha i},
\nabla\langle \nabla x_{\alpha}, \nabla u_i\rangle \rangle-
\int_{\om} x_{\alpha}\langle \nabla x_{\alpha}, \nabla
u_i\rangle\Delta u_i\\ \no &=&-\int_{\om} \langle\nabla h_{\alpha
i}, \nabla\langle \nabla x_{\alpha}, \nabla u_i\rangle\rangle - \fr
12\int_{\om} \langle \nabla x_{\alpha}^2, \nabla u_i\rangle\Delta
u_i, \en we have \be\no 1=\sum_{\alpha=1}^{n+1}\int_{\om}\langle
\nabla x_{\alpha}, \nabla
u_i\rangle^2=-\sum_{\alpha=1}^{n+1}\int_{\om} \langle\nabla
h_{\alpha i}, \nabla\langle \nabla x_{\alpha}, \nabla
u_i\rangle\rangle \en which gives \be\no
\sum_{\alpha=1}^{n+1}||\nabla\langle \nabla x_{\alpha}, \nabla
u_i\rangle||^2>0. \en It then follows from (2.77) that $
\laa_i^{1/(l-1)}-(n-2)>0$ and
\be\no 1&=& -\sum_{\alpha=1}^{n+1}\int_{\om} \langle\nabla h_{\alpha
i}, \nabla\langle \nabla x_{\alpha}, \nabla u_i\rangle\rangle\\
\no&\leq&\fr
12\sum_{\alpha=1}^{n+1}\left(\left(\laa_i^{1/(l-1)}-(n-2)\right)||\nabla
h_{\alpha i}||^2+\fr 1{\laa_i^{1/(l-1)}-(n-2)}||\nabla\langle \nabla
x_{\alpha}, \nabla u_i\rangle\rangle||^2\right)\\ \no &=& \fr 12+\fr
12\sum_{\alpha=1}^{n+1}\left((\laa_i^{1/(l-1)}-(n-2)\right)||\nabla
h_{\alpha i}||^2. \en Thus, we have \be
-\sum_{\alpha=1}^{n+1}||\nabla h_{\alpha i}||^2\leq -\fr
1{\laa_i^{1/(l-1)}-(n-2)} \en and consequently, one has \be
\sum_{\alpha=1}^{n+1}||W_{\alpha i}||^2=\sum_{\alpha=1}^{n+1}\left(
||x_{\alpha}\na u_i||^2-||\nabla h_{\alpha i}||^2\right)\leq 1 -\fr
1{\laa_i^{1/(l-1)}-(n-2)}. \en From $b_{\alpha ij}=b_{\alpha ji}, \
c_{\alpha ij}=-c_{\alpha ji}$, we have \be 2\sum_{i,
j=1}^k(\laa_{k+1}-\laa_i)^2b_{\alpha ij}c_{\alpha
ij}=-2\sum_{i,j=1}^k (\laa_{k+1}-\laa_i)(\laa_i-\laa_j)b_{\alpha
ij}c_{\alpha ij} \en \be \delta\sum_{i,
j=1}^k(\laa_{k+1}-\laa_i)^2(\laa_i-\laa_j)b_{\alpha
ij}^2=-\delta\sum_{i,j=1}^k(\laa_{k+1}-\laa_i)(\laa_i-\laa_j)^2b_{\alpha
ij}^2. \en
Combining  (2.54), (2.70) and (2.72), we get
\be& &
\sum_{\alpha=1}^{n+1}(\laa_{k+1}-\laa_i)||\nabla \phi_{\alpha i}||^2\\ \no &\leq&
\sum_{j=0}^{l-2}|a_j|\laa_i^{ j/(l-1)}
+\sum_{\alpha=1}^{n+1}\laa_i||W_{\alpha i}||^2+\sum_{\alpha=1}^{n+1}\sum_{j=1}^k(\laa_i-\laa_j)b_{\alpha ij}^2.
\en
We have by substituting (2.83) into (2.76)  that
\be & &
(\laa_{k+1}-\laa_i)^2\left(n+2\sum_{\alpha=1}^{n+1}\sum_{j=1}^k b_{\alpha ij}c_{\alpha ij}\right)
\\ \no &\leq &
\delta (\laa_{k+1}-\laa_i)^2\left( \sum_{j=0}^{l-2}|a_j|\laa_i^{ j/(l-1)}
+\sum_{\alpha=1}^{n+1}\sum_{j=1}^k(\laa_i-\laa_j)b_{\alpha ij}^2\right)
\\ \no & & +\fr{\laa_{k+1}-\laa_i}{\delta}\sum_{\alpha=1}^{n+1}\left(||Z_{\alpha i}||^2-\sum_{j=1}^kc_{\alpha ij}^2\right)
+\sum_{\alpha=1}^{n+1}(\laa_{k+1}-\laa_i)^2(\delta \laa_i+n-2) ||W_{\alpha i}||^2.
\en
Hence, by summing over $i$ from 1 to $k$ and noticing (2.78), (2.80), (2.81) and (2.82), we infer
\be \no & &
n\sum_{i=1}^k(\laa_{k+1}-\laa_i)^2
\\ \no &\leq & \delta (\laa_{k+1}-\laa_i)^2 \sum_{j=0}^{l-2}|a_j|\laa_i^{ j/(l-1)}+
\sum_{i=1}^k\fr{\laa_{k+1}-\laa_i}{\delta}\sum_{\alpha=1}^{n+1}||Z_{\alpha i}||^2\\ \no & &
 +\sum_{i=1}^k\sum_{\alpha=1}^{n+1}(\laa_{k+1}-\laa_i)^2(\delta \laa_i+n-2) ||W_{\alpha i}||^2
\\ \no &\leq &
\delta \sum_{i=1}^k(\laa_{k+1}-\laa_i)^2 \sum_{j=0}^{l-2}|a_j|\laa_i^{ j/(l-1)} + \sum_{i=1}^k\fr{\laa_{k+1}-\laa_i}{\delta}\left(\laa_i^{1/(l-1)}+\fr{(n-2)^2}4\right)
\\ \no & & +\sum_{i=1}^k(\laa_{k+1}-\laa_i)^2(\delta \laa_i+n-2)\left( 1
-\fr 1{\laa_i^{1/(l-1)}-(n-2)}\right).
\en
That is
\be  & &
\sum_{i=1}^k(\laa_{k+1}-\laa_i)^2\left(2+\fr{n-2}{\laa_i^{1/(l-1)}-(n-2)}\right)
\\ \no &\leq &
 \delta \sum_{i=1}^k(\laa_{k+1}-\laa_i)^2H_i+
\fr 1{\delta} \sum_{i=1}^k(\laa_{k+1}-\laa_i)\left(\laa_i^{1/(l-1)}+\fr{(n-2)^2}4\right),
\en
where $H_i$ is given by (1.22).
 This completes the proof of Theorem 1.2.

\vskip0.3cm
{\it Proof of Corollary 1.2.} Taking
\be\no & &
\delta =\fr{\left\{ \sum_{i=1}^k(\laa_{k+1}-\laa_i)\left(\laa_i^{1/(l-1)}+\fr{(n-2)^2}4\right)\right\}^{1/2}}{\left\{
 \sum_{i=1}^k(\laa_{k+1}-\laa_i)^2H_i\right\}^{1/2}},
\en
in  (1.21), we have
\be & &
\sum_{i=1}^k(\laa_{k+1}-\laa_i)^2\left(2+\fr{n-2}{\laa_i^{1/(l-1)}-(n-2)}\right)
\\ \no &\leq &
2\left\{ \sum_{i=1}^k(\laa_{k+1}-\laa_i)^2H_i \right\}^{1/2}
\times \left\{ \sum_{i=1}^k(\laa_{k+1}-\laa_i)\left(\laa_i^{1/(l-1)}+\fr{(n-2)^2}4\right)\right\}^{1/2},
\en
Since
\be\no
2+\fr{n-2}{\laa_i^{1/(l-1)}-(n-2)}\geq 2+\fr{n-2}{\laa_k^{1/(l-1)}-(n-2)}=S_k, \ i=1,\cdots, k,
\en
we have
\be
\sum_{i=1}^k(\laa_{k+1}-\laa_i)^2\left(2+\fr{n-2}{\laa_i^{1/(l-1)}-(n-2)}\right)\geq S_k\sum_{i=1}^k(\laa_{k+1}-\laa_i)^2
\en
and we infer from Lemma 2.1 that
\be
& & \left\{ \sum_{i=1}^k(\laa_{k+1}-\laa_i)^2H_i \right\}
\times \left\{ \sum_{i=1}^k(\laa_{k+1}-\laa_i)\left(\laa_i^{1/(l-1)}+\fr{(n-2)^2}4\right)\right\}
\\ \no &\leq& \left\{ \sum_{i=1}^k(\laa_{k+1}-\laa_i)^2 \right\}
\times \left\{ \sum_{i=1}^k(\laa_{k+1}-\laa_i)H_i\left(\laa_i^{1/(l-1)}+\fr{(n-2)^2}4\right)\right\}
\\ \no &=& \left\{ \sum_{i=1}^k(\laa_{k+1}-\laa_i)^2 \right\}
\times \left\{ \sum_{i=1}^k(\laa_{k+1}-\laa_i)T_i\right\},
\en
where $S_k$ and $T_i$ are defined as in (1.26).
Substituting (2.87) and (2.88) into (2.86), one gets
\be \no
\sum_{i=1}^k(\laa_{k+1}-\laa_i)^2\leq
\fr 4{S_k^2} \sum_{i=1}^k(\laa_{k+1}-\laa_i)T_i,
\en
where $A_{k+1}$ and $B_{k+1}$ are given by (1.25).
Solving this quadratic polynomial about $\laa_{k+1}$, we get (1.24).

\vskip0.3cm {\it Acknowledgements.} The fourth
author would like to thank the Max Planck Institute for Mathematics
in the Sciences for its hospitality and CAPES.

\vskip0.5cm
\noindent J\"urgen Jost ( jost@mis.mpg.de ), Xianqing Li-Jost (
xli-jost@mis.mpg.de )

\noindent Max Planck Institute for Mathematics in the Sciences,
04103 Leipzig, Germany

\vskip0.3cm \noindent Qiaoling Wang (
wang@mat.unb.br )

\noindent Departamento de Matem\'atica, UnB, 70910-900,
Bras\'{\i}lia-DF, Brazil

\vskip0.3cm
\noindent Changyu Xia ( xia@mat.unb.br ) \noindent Max Planck Institute for Mathematics in the Sciences,
04103 Leipzig, Germany,
and
Departamento de Matem\'atica, UnB, 70910-900, Bras\'{\i}lia-DF,
Brazil

\end{document}